\documentclass[12pt]{amsart}
\usepackage{amscd,amssymb}
\usepackage[graph,frame,poly,arc]{xy}  
\usepackage[plainpages,backref,urlcolor=blue]{hyperref}

\theoremstyle{plain}
\newtheorem{thm}[subsection]{Theorem}

\newtheorem{cor}[subsection]{Corollary}

\theoremstyle{definition}
\newtheorem{rk}[subsection]{Remark}

\newtheorem{ex}[subsection]{Example}
\newtheorem{conj}[subsection]{Conjecture}
\newtheorem{question}[subsection]{Question}

\numberwithin{equation}{section}
\newcommand{\C}{\mathbb{C}}
\newcommand{\PP}{\mathbb{P}}

\begin{document}

\title [Graded Betti numbers and total Tjurina numbers]
{Graded Betti numbers of the Jacobian algebra and total Tjurina numbers of plane curves}

\author[Alexandru Dimca]{Alexandru Dimca$^1$}
\address{Universit\'e C\^ ote d'Azur, CNRS, LJAD, France and Simion Stoilow Institute of Mathematics,
P.O. Box 1-764, RO-014700 Bucharest, Romania}
\email{Alexandru.Dimca@univ-cotedazur.fr}

\author[Gabriel Sticlaru]{Gabriel Sticlaru}
\address{Faculty of Mathematics and Informatics,
Ovidius University
Bd. Mamaia 124, 900527 Constanta,
Romania}
\email{gabriel.sticlaru@gmail.com }

\thanks{$^1$ partial support from the project ``Singularities and Applications'' - CF 132/31.07.2023 funded by the European Union - NextGenerationEU - through Romania's National Recovery and Resilience Plan.}



\begin{abstract}

In this paper we compute an explicit
closed formula for the total Tjurina number $\tau(C)$ of a reduced projective plane curve $C$ in terms of the graded Betti numbers of the corresponding Jacobian algebra.
This formula allows a completely new view point on the classical upper bounds for the total Tjurina number $\tau(C)$
of a plane curve $C$ given by A. du Plessis and C. T. C. Wall. This approach yields in particular a new necessary condition for a set of positive integers to be the graded Betti numbers of the Jacobian algebra of a reduced plane curve.
 
\end{abstract}

\maketitle

\section{Introduction} 

Let $S=\C[x,y,z]$ be the polynomial ring in three variables $x,y,z$ with complex coefficients, and let $C:f=0$ be a reduced curve of degree $d\geq 3$ in the complex projective plane $\PP^2$. 
We denote by $J_f$ the Jacobian ideal of $f$, i.e. the homogeneous ideal in $S$ spanned by the partial derivatives $f_x,f_y,f_z$ of $f$, and  by $M(f)=S/J_f$ the corresponding graded quotient ring, called the Jacobian (or Milnor) algebra of $f$.

The study of the Jacobian module $M(f)$ of a reduced plane curve
$C:f=0$ has become a central tool in understanding the geometry and
singularities of $C$.  
Among the numerical invariants extracted from $M(f)$, the \emph{type}
\[
t(C)=d_1+d_2+1-d \geq 0
\]
plays a particularly important role, where $d_1\le d_2\le\cdots\leq d_m$ are the
exponents of $C$, see for instance \cite{ADP,Type3} and $d=\deg f$. 
 We recall that $t(C)=0$ (resp. $t(C)=1$) if and only if the curve $C$ is free (resp. plus-one generated), see \cite{Abe,ADP}.
 After introducing  the type of a plane curve, T.~Abe,
P.~Pokora and the first author,  have obtained a classification into two classes
for curves with $t(C)=2$, see \cite{ADP}.
Later, we have established a classification into four
classes for curves with $t(C)=3$, see \cite{Type3}. 

In this paper we  show that in general the type $t$ curves fall naturally into $2^{t-1}$ classes, see Section 2. In Section 3, we  a closed formula for the corresponding Tjurina number in terms of the graded Betti numbers of the Jacobian algebra $M(f)$, see Theorem \ref{thm1}, which is the main result of this note. This formula is simplified by the use of the type $t=t(C)$.

As an example, we study these
classes associated to partitions in the case of reduced plane curves of type $t(C)=4$ in Section 4.
We show that exactly eight classes occur, 
and for each of them we show the non-emptiness by providing at least one example, see Example \ref{ex1}.

In Section 5 we look from a completely new view point at the classical upper bounds for the total Tjurina number $\tau(C)$
of a plane curve $C$ given by A. du Plessis and C.T.C. Wall in terms of $d$ and $d_1$. We also revisit the notions of nearly free curve, minimal plus-one generated curve and maximal Tjurina curve from this view point. The use of the type $t=t(C)$ of the curve $C$ makes the formulas much easier to state and also to use, see for instance the formulas in Corollary \ref{cor2}. 

This approach yields in particular a {\it new necessary condition for a set of positive integers to be the graded Betti numbers of the Jacobian algebra of a reduced plane curve}, see Remark \ref{rk2}.

By analogy to the characterization of nearly free curves recalled in 
Corollary \ref{cor3} (1), one may introduce the class of
{\it nearly maximal Tjurina curves} and describe the corresponding graded Betti numbers, see Corollary \ref{cor5}.

All the computations in the examples in this note were done using the software CoCoA \cite{CoCoA} and SINGULAR \cite{Singular}.

\section{Graded Betti numbers and partitions}

Consider  an ordered partition
\[
\pi  \  :  \  t=\epsilon_1+\cdots+\epsilon_s,
\]
where $t, \epsilon_1, \cdots ,\epsilon_s$ are strictly positive integers.
 
For integers $t\ge 1$ and $1\le s\le t$, we denote by
$
p(t,s)
$
the number of such ordered partitions of $t$ with exactly $s$ parts.  
Since such partitions correspond bijectively to choices of $s-1$ points
among the $t-1$ points $\{1,2, \ldots,t-1\}$ on the real line,
 one has
\[
p(t,s)=\binom{t-1}{s-1}.
\]
The total number of partitions of $t$ is therefore
\[
p(t)=\sum_{s=1}^{t} p(t,s)
     =\sum_{s=1}^{t} \binom{t-1}{s-1}
     =2^{t-1}.
\]
Consider now the general form of the minimal resolution of the Milnor algebra $M(f)$
of a reduced plane curve $C:f=0$ which is assumed to be not free:
\begin{equation}
\label{res2A}
0 \longrightarrow \bigoplus_{j=1}^{m-2} S(-c_j)
   \longrightarrow \bigoplus_{i=1}^{m} S(1-d-d_i)
   \longrightarrow S^3(1-d)
   \longrightarrow S,
\end{equation}
with $c_1\leq ...\leq c_{m-2}$ and $d_1\leq ...\leq d_m$.
It follows from \cite[Lemma 1.1]{HS} that
\begin{equation}
\label{res2B}
c_j = d + d_{j+2} - 1 + \epsilon_j,
\qquad j = 1,\dots,m-2,
\end{equation}
for some integers $\epsilon_j \ge 1$. 

We call the ordered sequence of degrees $d_i$ for $i=1, \ldots,m$ and the integers $\epsilon_j$  for $j=1, \ldots, m-2$ for $m \geq 3$ the {\it graded Betti numbers of the Jacobian algebra} $M(f)$, since they determine and are determined by the usual
graded Betti numbers of the Jacobian algebra $M(f)$ as defined for instance in \cite{Eis}.

Using \cite[Formula (13)]{HS}, one obtains the relation
\begin{equation}
\label{res2C}
d_1 + d_2 = d - 1 + \sum_{j=1}^{m-2} \epsilon_j.
\end{equation}
Hence we get an ordered partition
\begin{equation}
\label{res2C1}
 \pi _C \  :  \ t(C) =  \sum_{j=1}^{m-2} \epsilon_j.
\end{equation}
In other words, each curve $C$ with $t(C)=t$ induces an ordered partition $\pi_C$ of the integer $t\geq 1$. It is natural therefore to divide the curves of type $t$ into $2^{t-1}$ classes, each class, say denoted by $tT_{\pi}$, being determined by an ordered partition $\pi$ on $t$. More precisely,
we say that $C$ has type $tT_{\pi}$ if $t(C)=t$ and $\pi_C=\pi$.

We put forth the following.

\begin{conj}
\label{conj1}
 For each $t\geq 1$ and each ordered partition $\pi$ of $t$, the set of
 plane curves of class $tT_{\pi}$ is not empty.
\end{conj}
This conjecture holds for $t=1$ since in this case we have only one class and for $t=2$ (resp. $t=3$) as shown in
\cite{ADP} (resp. \cite{Type3}). Example \ref{ex1} below shows that this conjecture also holds for $t=4$. 
A much more difficult problem seems to be the following.

\begin{question}
\label{q1}
Find the necessary and sufficient conditions on the ordered sequence of degrees $d_i$ for $i=1, \ldots,m$ and on the integers $\epsilon_j\geq 1$  for $j=1, \ldots, m-2$, where $m \geq 3$, such that there is a reduced, non free curve $C:f=0$ with these data the graded Betti numbers of the Jacobian algebra $M(f)$.
\end{question}
Note that the degree of such a curve $C$ is determined by \eqref{res2C} and the inequality $c_j \leq c_{j+1}$  stated in \eqref{res2A} yields via
\eqref{res2B} the necessary conditions
\begin{equation}
\label{ineq}
d_{j+2}+\epsilon_j \leq d_{j+3}+\epsilon_{j+1}
\end{equation}
for $j=1, \ldots, m-3$. We also know that
\begin{equation}
\label{boundG}
d_1\leq d_2 \leq d_3 \leq d-1,
\end{equation}
see for instance \cite[Theorem 2.4]{3syz} and
\begin{equation}
\label{boundG2}
 d_m \leq 2d-4,
\end{equation}
see for instance \cite[inequality (2.6)]{CMreg}.
Let us consider the following expression involving the graded Betti numbers
\begin{equation}
\label{boundG31}
E=d_1^2+td_2-\frac{t^2+ \sum_{j=1}^n\epsilon_j^2}{2}-\sum_{j=1}^n\epsilon_jd_{j+2},
\end{equation}
where $t=t(C)$ as in \eqref{res2C1}.
In this note we prove the following inequality
\begin{equation}
\label{boundG3}
E\geq 0
\end{equation}
with equality if and only if $m=3$ and $d_2=d_3=d-1$,
which is a {\it new necessary condition to be satisfied by the graded Betti numbers of a reduced plane curve}. In other words, \eqref{boundG3} is not a consequence of the inequalities
\eqref{ineq}, \eqref{boundG} and \eqref{boundG2}, see for details 
Remark \ref{rk2} below.

%%%%%%%%%%%%%%%%%%%%%%%%%%%%%%%%%%%%%%%%%%%%%%%%%%%%
\section{Tjurina numbers for curves via graded Betti numbers of their Jacobian algebra}
%%%%%%%%%%%%%%%%%%%%%%%%%%%%%%%%%%%%%%%%%%%%%%%%%%%%

\begin{thm}
\label{thm1}
Let $C:f=0$ be a reduced plane curve of degree $d$ and type $t(C)=t\geq 1$ with exponents
$(d_{1},d_{2},\ldots,d_{m})$.
Let $\pi_C=(\epsilon_1,\dots,\epsilon_n)$ be the corresponding ordered partition of $t$, with $n=m-2$. Then one has
$$\tau(C)=d_1^2+d_1d_2+ d_2^2-t(d_1+d_2)-\sum_{i=3}^m\epsilon_{i-2}d_i+R,$$
where
$$R=\frac{t^2-\sum_{i=3}^m\epsilon_{i-2}^2}{2}= \sum_{3\leq i <j \leq m}
\epsilon_{i-2}\epsilon_{j-2} \geq 0.$$
In particular, $R=0$ if and only if $m=3$.
\end{thm}

\proof
By taking the homogeneous component of degree $k$ and using exactness in the resolution \eqref{res2A} we get the following equality, for $k \gg 0$.
$$
\dim M(f)_k
=\binom{k+2}{2}
-3\binom{k+3-d}{2}
+\sum_{i=1}^{n+2}\binom{k+3-d-d_i}{2}-$$
$$
-\sum_{j=1}^{n}\binom{k+3-d-d_{j+2}-\epsilon_j}{2}.
$$
This is a polynomial in $k$ of degree at most $2$.  
The constant term of this polynomial is $\tau(C)$, and a direct computation, done using $d=d_1+d_2+1-t$, the equality \eqref{res2C1}, and grouping the terms in the above sum as follows
$$
\dim M(f)_k
=\binom{k+2}{2}
-3\binom{k+3-d}{2}
+\sum_{i=1}^{2}\binom{k+3-d-d_i}{2}+$$
$$+\sum_{j=1}^n\left( \binom{k+3-d-d_{j+2}}{2}-
\binom{k+3-d-d_{j+2}-\epsilon_j}{2}\right)
$$
yields the result.

\endproof

\begin{rk}
\label{rk1}
(i) Theorem \ref{thm1} is a generalization of \cite[Proposition 1.11]{ADP} which treats the case $t=2$ and of \cite[Theorem 3.1]{Type3} which treats the case $t=3$.

\medskip

\noindent (ii) When $C$ is a free curve of degree $d$ with exponents $(d_1,d_2)$, then
$d_1+d_2=d-1$ and 
$$\tau(C)=(d-1)^2-d_1(d-d_1-1)=(d_1+d_2)^2-d_1d_2=d_1^2+d_1d_2+ d_2^2,$$
see for instance Theorem \ref{thmCTC} below. Since in this case $t(C)=0$, we see that the equality in Theorem \ref{thm1} holds in this case as well with the obvious conventions
$$\sum_{i=3}^m\epsilon_{i-2}d_i=\sum_{i=3}^m\epsilon_{i-2}^2=0.$$
Indeed, in this case there are no $\epsilon_j$ since $n=m-2=0$.

\end{rk}

\begin{ex}
\label{ex0}
(i) For the partition $\pi=(t)$ of $t(C)=t$, we get the formula
$$\tau(C)=d_1^2+d_1d_2+d_2^2-t(d_1+d_2+d_3).$$
Indeed, in this case $m=3$ and $\epsilon_1=t$.
To have a geometric example,
if $C$ is a smooth curve, then the exponents of $C$ are
$$(d-1,d-1,d-1)$$
and the type of $C$ is 
$t(C)=d-1$. If we consider the partition $\pi=(d-1)$, then clearly the smooth curve $C$ has class $(d-1)T_{\pi}$.

\medskip

\noindent (ii) For the partition $\pi'=(1,1, \ldots,1)$ of $t(C)=t$, we get the formula
$$\tau(C)=d_1^2+d_1d_2+d_2^2-t(d_1+d_2)-\sum_{j=3}^{t+2}d_j+\frac{t(t-1)}{2}.$$
In this case $m=t+2$ and $\epsilon_1= \ldots = \epsilon_t=1$.
Note that the maximal Tjurina curves discussed in Corollary \ref{cor4} below are in this class of curves.

\end{ex}

%%%%%%%%%%%%%%%%%%%%%%%%%%%%%%%%%%%%%%%%%%%%%%
\section{Examples for curves of type $t(C)=4$}
%%%%%%%%%%%%%%%%%%%%%%%%%%%%%%%%%%%%%%%%%%%%%%%%%%%%
The following is a direct consequence of Theorem \ref{thm1}.
Consider all the possible partitions of $t=4$, namely:
$$\pi_1=(4), \ \pi_2=(3,1), \  \pi_3=(2,2), \ \pi_4=(1,3), \ \pi_5=(2,1,1), \ \pi_6=(1,2,1), $$
$$ \pi_7=(1,1,2) \text{ and } \pi_8=(1,1,1,1).$$
\begin{cor}
\label{cor1}
Let $C$ be a reduced plane curve of type $t(C)=4$ with exponents
$(d_{1},d_{2},\ldots,d_{m})$.  
Then $d_{1}+d_{2}=d+3$, and $C$ is in exactly one of the following eight classes.

\begin{enumerate}

\item[\textbf{Type $4T_{\pi_1}$.}]
$C$ is a $3$-syzygy curve and $\epsilon_{1}=4$.  
If $(d_{1},d_{2},d_{3})$ are the exponents of $C$, then
\[
\tau(C)=d_{1}^{2}+d_{1}d_{2}+d_{2}^{2}-4(d_{1}+d_{2}+d_{3}).
\]

\item[\textbf{Type $4T_{\pi_2}$.}]
$C$ is a $4$-syzygy curve and $\epsilon_{1}=3$, $\epsilon_{2}=1$.  
If $(d_{1},d_{2},d_{3},d_{4})$ are the exponents of $C$, then
$d_3<d_4-1$ and
\[
\tau(C)
= d_{1}^{2}+d_{1}d_{2}+d_{2}^{2}
-4d_{1}-4d_{2}
+3
-3d_{3}-d_{4}.
\]

\item[\textbf{Type $4T_{\pi_3}$.}]
$C$ is a $4$-syzygy curve and $\epsilon_{1}=2$, $\epsilon_{2}=2$.  
If $(d_{1},d_{2},d_{3},d_{4})$ are the exponents of $C$, then
\[
\tau(C)
= d_{1}^{2}+d_{1}d_{2}+d_{2}^{2}
-4d_{1}-4d_{2}
+4
-2d_{3}-2d_4.
\]

\item[\textbf{Type $4T_{\pi_4}$.}]
$C$ is a $4$-syzygy curve and $\epsilon_{1}=1$, $\epsilon_{2}=3$.  
If $(d_{1},d_{2},d_{3},d_{4})$ are the exponents of $C$, then
\[
\tau(C)
= d_{1}^{2}+d_{1}d_{2}+d_{2}^{2}
-4d_{1}-4d_{2}
+3
-d_{3}-3d_{4}.
\]

\item[\textbf{Type $4T_{\pi_5}$.}]
$C$ is a $5$-syzygy curve and $\epsilon_{1}=2$, $\epsilon_{2}=1$, $\epsilon_{3}=1$.  
If $(d_{1},d_{2},d_{3},d_{4},d_{5})$ are the exponents of $C$, then
$d_3<d_4$ and
\[
\tau(C)
= d_{1}^{2}+d_{1}d_{2}+d_{2}^{2}
-4d_{1}-4d_{2}
+5
-2d_{3}-d_{4}-d_{5}.
\]

\item[\textbf{Type $4T_{\pi_6}$.}]
$C$ is a $5$-syzygy curve and $\epsilon_{1}=1$, $\epsilon_{2}=2$, $\epsilon_{3}=1$.  
If $(d_{1},d_{2},d_{3},d_{4},d_{5})$ are the exponents of $C$, then
$d_4<d_5$ and
\[
\tau(C)
= d_{1}^{2}+d_{1}d_{2}+d_{2}^{2}
-4d_{1}-4d_{2}
+5
-d_{3}-2d_{4}-d_{5}.
\]

\item[\textbf{Type $4T_{\pi_7}$.}]
$C$ is a $5$-syzygy curve and $\epsilon_{1}=1$, $\epsilon_{2}=1$, $\epsilon_{3}=2$.  
If $(d_{1},d_{2},d_{3},d_{4},d_{5})$ are the exponents of $C$, then
\[
\tau(C)
= d_{1}^{2}+d_{1}d_{2}+d_{2}^{2}
-4d_{1}-4d_{2}
+5
-d_{3}-d_{4}-2d_{5}.
\]

\item[\textbf{Type $4T_{\pi_8}$.}]
$C$ is a $6$-syzygy curve and $\epsilon_{1}=\epsilon_{2}=\epsilon_{3}=\epsilon_{4}=1$.  
If $(d_{1},d_{2},d_{3},d_{4},d_{5},d_{6})$ are the exponents of $C$, then
\[
\tau(C)
= d_{1}^{2}+d_{1}d_{2}+d_{2}^{2}
-4d_{1}-4d_{2}
+6
-(d_{3}+d_{4}+d_{5}+d_{6}).
\]

\end{enumerate}
\end{cor}
\proof
The only points that need additional explanations are the following.
The inequality $d_3 <d_4-1$ in Type $4T_{\pi_2}$ and the inequality $d_3 <d_4$ in Type $4T_{\pi_5}$ come from the inequality stated in
\eqref{ineq} for $j=1$
The inequality $d_4 <d_5$ in Type $4T_{\pi_6}$ comes from the similar inequality stated in
\eqref{ineq} for $j=2$.
\endproof

%%%%%%%%%%%%%%%%%%%%%%%%%%%%%%%%%%%%%%%%%%%%%%%%%%%%
%\subsection*{Examples of curves of type $t(C)=4$}
%%%%%%%%%%%%%%%%%%%%%%%%%%%%%%%%%%%%%%%%%%%%%%%%%%%%

We give next explicit examples of reduced plane curves of
type $t(C)=4$, one for each of the eight possible classes listed above. For each example we list the
degree $d$, the Tjurina number $\tau(C)$, the exponents of $C$, and the
minimal free resolution of the Milnor algebra $M(f)$.

\begin{ex}
\label{ex1}

\noindent\textbf{(i) Type $4T_{\pi_1}$}
Consider the curve
\[
C:f=x^{5}+y^{5}+z^{5}=0.
\]
Then $d=5$, $t(C)=4$, $m=3$, $\tau(C)=0$, and the exponents are
$(4,4,4)$.
The minimal free resolution of $M(f)$ is
\[
0 \to S(-12) \to S(-8)^{3} \to S(-4)^{3} \to S.
\]

\medskip
\noindent\textbf{(ii) Type $4T_{\pi_2}$}
Consider the curve
\[
C:f=x^{5}+xyz^{3}+z^{5}=0.
\]
Then $d=5$, $t(C)=4$, $m=4$, $\tau(C)=1$, and the exponents are
$(4,4,4,6)$.
The minimal free resolution of $M(f)$ is
\[
0 \to S(-11)^{2} \to S(-8)^{3} \oplus S(-10) \to S(-4)^{3} \to S.
\]

\medskip
\noindent\textbf{(iii) Type $4T_{\pi_3}$}

Consider the curve
\[
C:f=(x^{3}+y^{3})^{3}+(y^{3}+z^{3})^{3}=0.
\]
Then $d=9$, $t(C)=4$, $m=4$, $\tau(C)=36$, and the exponents are
$(4,8,8,8)$.
The minimal free resolution of $M(f)$ is
\[
0 \to S(-18)^{2} \to S(-12)\oplus S(-16)^{3} \to S(-8)^{3} \to S.
\]

\medskip
\noindent\textbf{(iv) Type $4T_{\pi_4}$}

Consider the curve
\[
C:f=x^{6}+y^{4}z^{2}+z^{4}x^{2}=0.
\]
Then $d=6$, $t(C)=4$, $m=4$, $\tau(C)=8$, and the exponents are
$(4,5,5,5)$.
The minimal free resolution of $M(f)$ is
\[
0 \to S(-11)\oplus S(-13) \to S(-9)\oplus S(-10)^{3} \to S(-5)^{3} \to S.
\]

\medskip
\noindent\textbf{(v) Type $4T_{\pi_5}$}
Consider the curve
\[
C:f=x^{7}+y^{7}-4x^{2}y^{2}z^{3}=0.
\]
Then $d=7$, $t(C)=4$, $m=5$, $\tau(C)=14$, and the exponents are
$(5,5,6,7,7)$.
The minimal free resolution of $M(f)$ is
\[
0 \to S(-14)^{3} \to S(-11)^{2}\oplus S(-12)\oplus S(-13)^{2}
\to S(-6)^{3} \to S.
\]

\medskip
\noindent\textbf{(vi) Type $4T_{\pi_6}$}
Consider the curve
\[
C:f=z(x+y-z)(2x+3y-5z)(3x+5y-7z)(x^{3}+y^{3}-z^{3})=0.
\]
Then $d=7$, $t(C)=4$, $m=5$, $\tau(C)=19$, and the exponents are
$(5,5,5,5,6)$.
The minimal free resolution of $M(f)$ is
\[
0 \to S(-12)\oplus S(-13)^{2} \to S(-11)^{4}\oplus S(-12)
\to S(-6)^{3} \to S.
\]

\medskip
\noindent\textbf{(vii) Type $4T_{\pi_7}$}
Consider the curve
\[
C:f=(x^{4}+y^{4})^{2}-4x^{2}y^{2}z^{4}=0.
\]
Then $d=8$, $t(C)=4$, $m=5$, $\tau(C)=28$, and the exponents are
$(5,6,6,6,6)$.
The minimal free resolution of $M(f)$ is
\[
0 \to S(-14)^{2}\oplus S(-15) \to S(-12)\oplus S(-13)^{4}
\to S(-7)^{3} \to S.
\]

\medskip
\noindent\textbf{(viii) Type $4T_{\pi_8}$}
Consider the curve
\[
C:f=x^{3}y^{2}+y^{3}z^{2}+x^{2}z^{3}=0.
\]
Then $d=5$, $t(C)=4$, $m=6$, $\tau(C)=6$, and the exponents are
$(4,4,4,4,4,4)$.
The minimal free resolution of $M(f)$ is
\[
0 \to S(-9)^{4} \to S(-8)^{6} \to S(-4)^{3} \to S.
\]
\end{ex}

\section{On a result by A. du Plessis and C.T.C. Wall}

We recall first the following result due to du Plessis and Wall, see \cite[Theorem 3.2]{duPCTC} as well as \cite{E} for an alternative approach.
\begin{thm}
\label{thmCTC}
For positive integers $d$ and $d_1$, define two new integers by 
$$\tau(d,d_1)_{min}=(d-1)(d-d_1-1)  \text{ and } 
\tau(d,d_1)_{max}= (d-1)^2-d_1(d-d_1-1).$$ 
If $C:f=0$ is a reduced curve of degree $d$ in $\PP^2$ and  $d_1$ is the first exponent of $f$,  then
$$\tau(d,d_1)_{min} \leq \tau(C) \leq \tau(d,d_1)_{max}$$
and the equality $\tau(C) = \tau(d,d_1)_{max}$ holds if and only if
$2d_1<d$ and $C$ is a free curve.
Moreover, for $2d_1 \geq d$, the stronger inequality
$\tau(C) \leq \tau'(d,d_1)_{max}$ holds, where
$$\tau(d,d_1)'_{max}=\tau(d,d_1)_{max} - \binom{2d_1+2-d}{2}.$$
\end{thm}

We recall that a {\it maximal Tjurina curve } $C$ is a curve such that
$2d_1 \geq d$ and
$\tau(C)=\tau (d,d_1)'_{max}$, 
see \cite{maxTjurina}. As an example, the curve in Example \ref{ex1} (viii) is a maximal Tjurina curve with $d=5$ and $d_1=4$.
Theorem \ref{thm1} yields a precise way to measure the difference between $\tau(C)$ and $ \tau (d,d_1)_{max}$ (resp. $\tau (d,d_1)'_{max}$), and in this way gives, in particular, a proof of the claims involving 
these invariants in Theorem \ref{thmCTC}.

\begin{cor}
\label{cor2}
Let $C:f=0$ be a reduced plane curve of degree $d$ and type $t(C)=t\geq 1$ with exponents
$(d_{1},d_{2},\ldots,d_{m})$ 
and $\pi_C=(\epsilon_1,\dots,\epsilon_n)$  the corresponding ordered partition of $t$, with $n=m-2$. Then one has the following.

\begin{enumerate}

   \item For any curve $C$ as above, one has
$$ \tau (d,d_1)_{max}-\tau(C)=\frac{t^2+t}{2}+\frac{\sum_{j=1}^n(\epsilon_j^2-\epsilon_j)}{2}+\sum_{j=1}^n\epsilon_j(d_{j+2}-d_2).$$
    
    \item If in addition $2d_1 \geq d$, then
$$ \tau (d,d_1)'_{max}-\tau(C)=\frac{t^2-(t-(d_2-d_1))^2}{2}+\frac{\sum_{j=1}^n(\epsilon_j^2-\epsilon_j)+(d_2-d_1)}{2}+\sum_{j=1}^n\epsilon_j(d_{j+2}-d_2).$$
\end{enumerate}
Moreover, all the terms in the sums in right hand sides in (1) and (2) are positive.

\end{cor}
\proof
We give the proof only in the case (2), the case (1) being similar and simpler.
Since $d_1+d_2 \geq 2d_1 >d-1$, it follows that 
$t=t(C)\geq 1$. Using the formula for $\tau(C)$ given in Theorem \ref{thm1}, the equality \eqref{res2C1} and replacing $2d_1+2-d$ in the binomial coefficient by
$t-(d_2-d_1)+1$, we get
$$(d-1)^2-d_1(d-d_1-1)-\binom{t-(d_2-d_1)+1}{2}-\tau(C)=$$
$$=\frac{t^2-(t-(d_2-d_1))^2}{2}+\frac{1}{2}\sum_{j=1}^n(\epsilon_j^2-\epsilon_j)+\frac{d_2-d_1}{2}+\sum_{j=1}^n\epsilon_j(d_{j+2}-d_2).$$
The right hand side is a sum where all the terms are positive. Indeed,
note that
$$t^2-(t-(d_2-d_1)^2=(d_2-d_1)(2t -(d_2-d_1)=2(d_2-d_1)(2d_1+1-d)\geq0.$$
Hence this sum is equal to 0 if and only if all the terms are equal to 0.
This implies that
$$d_1=d_2= \ldots =d_m$$
and $\epsilon_j=1$ for all $j's$. Finally, the equality \eqref{res2C} yields
$m=2d_1-d+3$, which completes our proof.
\endproof
\begin{rk}
\label{rk2}
A direct computation using the formula for $\tau(C)$ given in Theorem \ref{thm1} yields
$$E:=\tau(C)-\tau(d,d_1)_{min}=d_1^2+td_2-\frac{t^2+ \sum_{j=1}^n\epsilon_j^2}{2}-\sum_{j=1}^n\epsilon_jd_{j+2}.$$
For $m=3$ the inequality $E \geq 0$ holds for obvious reasons, since one has $\epsilon_1=t$
and hence
$$E=d_1(d_1-t)+t(d_1+d_2-t-d_3)=d_1(d_1-t)+t(d-1-d_3) \geq 0.$$
Indeed, $d_1-t=d-1-d_2 \geq 0$ and $d-1-d_3 \geq 0$ by \eqref{boundG}. In particular, the equality $E=0$ is equivalent to 
$$d_2=d_3=d-1$$
when $m=3$.
On the other hand, for $m \geq 4$, the inequality $E \geq 0$
coming from Theorem \ref{thmCTC} is a new necessary condition to be satisfied by the exponents
$(d_{1},d_{2},\ldots,d_{m})$ 
and $(\epsilon_1,\dots,\epsilon_n)$, the corresponding ordered partition of $t$. In fact using \cite[Theorem 3.5 (1)]{3syz} we get the following more precise inequality for the graded Betti numbers of any reduced curve $C$
\begin{equation}
\label{ineqMIN}
E>0
\end{equation}
unless $m=3$ and $d_2=d_3=d-1$ when $E=0$.
To show that this inequality is not a consequence of the inequalities
\eqref{ineq}, \eqref{boundG} and \eqref{boundG2}, one may consider
 the data $d_1=3$, $d_2=d_3=8$, $d_4=15$,
$\epsilon_1=\epsilon_2=1$ which  verify the inequalities
\eqref{ineq}, \eqref{boundG} and \eqref{boundG2}, but not the
above inequality. 
 Indeed, for these data one gets $E=-1 <0$, which shows that there is no reduced curve $C$ with these Betti numbers.
\end{rk}

Recall that an $3$-syzygy curve $C$ is called nearly free (resp. minimal plus-one generated) if $d_1+d_2=d$  and $d_2=d_3$ (resp. 
$d_1+d_2=d$ and $d_2+1=d_3$, see \cite{DStRIMS,Brian}.
Using Corollary \ref{cor2} (1) we can give new proofs for the following known results. Indeed, the claim (1) below occurs in \cite[Theorem 1.3]{Dmax}, and the claim (2) below occurs in \cite[Theorem 1.5]{Brian}.

\begin{cor}
\label{cor3}
Let $C:f=0$ be a reduced plane curve of degree $d$ with exponents
$(d_{1},d_{2},\ldots,d_{m})$. Then one has the following.

\begin{enumerate}

   \item $C$ is a nearly free curve if and only if 
$ \tau(C)=\tau (d,d_1)_{max}-1$.
    
    \item $C$ is a  minimal plus-one generated curve if and only if 
$ \tau(C)=\tau (d,d_1)_{max}-2$.
    
\end{enumerate}
\end{cor}
\proof
If the curve $C$ is nearly free or respectively minimal plus-one generated, it follows from the general formula for the Tjurina number of
a 3-syzygy curve given in \cite{3syz} that $\tau(C)$ satisfies the claimed relations.

Conversely, if $\tau(C)$ satisfies the claimed relations,
it follows from Theorem \ref{thmCTC} that the curve $C$ is not free, and hence $t(C) \geq 1$.

It is clear that
$$\tau (d,d_1)_{max}- \tau(C)=1$$
implies $t=1$, and hence $C$ is a plus-one generated curve, in particular $m=3$. Then we get $\epsilon_1=1$ and $d_2=d_3$, which tells us that $C$ is a nearly free curve.

Next, if we have
$$\tau (d,d_1)_{max}- \tau(C)=2,$$
then again $t=1$ and $m=3$ as above. Then we get either
$\epsilon_1=1$ or $\epsilon_1=2$. In the first case we get
$d_3=d_2+1$, that is $C$ is a  minimal plus-one generated curve.
In the second case, we get $d_2=d_3$, that is $C$ is a nearly free curve. But this is impossible, since it contradicts (1).

\endproof

Using Corollary \ref{cor2} (2) we can give
a completely new proof of \cite[Theorem 3.1]{maxTjurina}, that characterizes the maximal Tjurina curves in terms of their resolutions.

\begin{cor}
\label{cor4}
Let $C:f=0$ be a reduced plane curve of degree $d$ and type $t(C)=t\geq 1$ with exponents
$(d_{1},d_{2},\ldots,d_{m})$ 
and $\pi_C=(\epsilon_1,\dots,\epsilon_n)$  the corresponding ordered partition of $t$, with $n=m-2$. 
 Then $C:f=0$ is a maximal Tjurina curve  if and only if  the exponents of $C$ satisfy
$$d_1=d_2= \ldots =d_m,$$
with $m=2d_1-d+3$ and 
$$\epsilon_1= \ldots = \epsilon_n=1.$$ 
\end{cor}
\proof
The sum in Corollary \ref{cor2} (2)  is equal to 0 if and only if all the terms are equal to 0.
This implies that
$$d_1=d_2= \ldots =d_m$$
and $\epsilon_j=1$ for all $j's$. Finally, the equality \eqref{res2C} yields
$m=2d_1-d+3$, which completes our proof.
\endproof

We have the following analog of Corollary \ref{cor3} (1), where $\tau (d,d_1)_{max}$ is replaced by $\tau (d,d_1)'_{max}$. We say that  $C:f=0$ is a {\it nearly maximal Tjurina curve} if $2d_1 \geq d$ and
$$\tau(C)=\tau (d,d_1)'_{max}-1.$$

\begin{cor}
\label{cor5}
Let $C:f=0$ be a reduced plane curve of degree $d$ and type $t(C)=t\geq 1$ with exponents
$(d_{1},d_{2},\ldots,d_{m})$ 
and $\pi_C=(\epsilon_1,\dots,\epsilon_n)$  the corresponding ordered partition of $t$, where $n=m-2$. Then  $C$ is a nearly maximal Tjurina curve if and only if  
one of  the following two situations occurs.

\begin{enumerate}

   \item $d_1=d_2=\ldots =d_{m-1}$, $d_m=d_1+1$ and
    $\epsilon_1= \ldots = \epsilon_n=1$, with $m=2d_1+3-d$.

   \item $d_1=d_2=\ldots =d_m$, $\epsilon_1= \ldots = \epsilon_{n-1}=1$
  and $\epsilon_n=2$, with $m=2d_1+2-d$.

\end{enumerate}

\end{cor}

\proof
Note that $t(C) \geq 1$ implies $m \geq 3$, which yields $2d_1 \geq d$ in both cases (1) and (2). If one of the situations (1) or (2) occurs, then
using Corollary \ref{cor2} (2) it follows that $C$ is a nearly maximal Tjurina curve.

Conversely, assume now that $C$ is a nearly maximal Tjurina curve.
Assume first that $d_1=d_2$. Then either all the $\epsilon _j$ are $1$,
and this leads to the case (1), or all the $\epsilon _j$ are $1$ except one, which is $2$. Then we have $d_1=d_2=\ldots =d_m$ and we must have
$\epsilon_n=2$ due to the inequality $c_{n-1} \leq c_n$ in \eqref{res2A}.
Hence this case corresponds to the claim (2) in our result.

Assume now that $d_1<d_2$. Then 
$$t^2-(t-(d_2-d_1)^2=2(d_2-d_1)(2d_1+1-d)\geq 2,$$
as we have seen in the proof of Corollary \ref{cor2}.
It follows that
$$\frac{t^2-(t-(d_2-d_1))^2}{2}+\frac{d_2-d_1}{2} >1$$
which is a contradiction. Hence this case cannot occur.
\endproof

\begin{rk}
\label{rk5}
To state the conditions on the partition $\pi_C=(\epsilon_1,\dots,\epsilon_n)$ that occurs in Corollary \ref{cor5}, it is convenient to introduce the invariant of $C$ given by
$$\Delta m=t-(m-2)=d_1+d_2-d-m+3,$$
which was already considered in \cite{SameDeg} when $d_1=d_2$.
Using \eqref{res2C1}, it follows that
$$\Delta m=\sum_{j=1}^{m-2}(\epsilon_j-1)\geq 0.$$
It follows that the condition on the partition $\pi_C$ in Corollary \ref{cor5} (1) is equivalent to
$$\Delta m =0,$$
and  the condition on the partition $\pi_C$ in Corollary \ref{cor5} (2) is equivalent to
$$\Delta m=1,$$
when the additional condition $d_1=d_2=\ldots =d_m$ holds.
With the terminology introduced in \cite{SameDeg} a curve $C$ satisfying $\pi_C$ in Corollary \ref{cor5} (2) is exactly a curve $C$ of type $(d,d_1,m)$ such that $\Delta m=1$.

\end{rk}

\begin{ex}
\label{ex5}
Here are some examples of curves satisfying the conditions in Corollary \ref{cor5}. 

(i) Consider the curve
$$C:f=xyz(x-z)(x-2z)(x-3z)(y-z)(y-2z)(y-3z)(x+y)(x+y-2z)=0$$
from \cite[Remark 4.3]{SameDeg}.
Then $C$ is a 4-syzygy curve with exponents $(6,6,6,7)$ and degree $d=11$. It follows that
$$\Delta m=d_1+d_2-d-m+3=6+6-11-4+3=0,$$
therefore by our Remark \ref{rk5} this curve satisfies the conditions in Corollary \ref{cor5} (1).

Similarly, the curve 
$$C':f'= xyz(x+y)(y+z)(z+x)(-x+y+z)(x-y+z)(x+y-z)=0$$
has exponents $(5,5,5,6)$ and $d=9$, hence again 
$$\Delta m=d_1+d_2-d-m+3=5+5-9-4+3=0.$$

The irreducible quartic 
$$C'': f''=x^4+x^3z+y^2z^2=0$$
has exponents $(2,2,3)$ and hence it is minimal plus-one generated as in Corollary \ref{cor3}. On the other hand it verifies
$$\Delta m==d_1+d_2-d-m+3=2+2-4-3+3=0.$$

(ii) Using the final part of Remark \ref{rk5} we have the following.
The curves $C$ and $C'$ in \cite[Example 4.5]{SameDeg} satisfy the conditions in Corollary \ref{cor5} (2) for $d=7$, $d_1=4$ and $m=3$.
The curves $C$  in \cite[Example 4.6]{SameDeg}  satisfies the conditions in Corollary \ref{cor5} (2) for $d=9$, $d_1=5$ and $m=3$.
The curves $C$  in \cite[Example 4.7]{SameDeg}  satisfies the conditions in Corollary \ref{cor5} (2) for $d=11$, $d_1=6$ and $m=3$.
To get examples of curves satisfying the conditions in Corollary \ref{cor5} (2) with $m\geq 4$, one may consider
the curves $C$, $C'$, $C'''$ and $C_6$ in \cite[Example 5.1]{SameDeg}, the curves $C_7$ and $C'$ in \cite[Example 5.2]{SameDeg}, the
curves $C$ and $C_8$ in \cite[Example 5.3]{SameDeg} and the curves $C'$ and  $C'''$ in \cite[Example 5.4]{SameDeg}.
\end{ex}

\section*{Conflict of Interests}
We declare that there is no conflict of interest regarding the publication of this paper.
\section*{Data Availability Statement}
We do not analyze or generate any datasets, because this work proceeds within a theoretical and mathematical approach.

\end{document}